\documentclass{amsart}

\newtheorem{theorem}{Theorem}[section]
\newtheorem{lemma}[theorem]{Lemma}
\newtheorem{corollary}[theorem]{Corollary}

\theoremstyle{definition}
\newtheorem{definition}[theorem]{Definition}

\theoremstyle{remark}
\newtheorem{remark}[theorem]{Remark}

\numberwithin{equation}{section}

\begin{document}

\title[A characterization of subspaces and quotients...]
{A characterization of subspaces and quotients of reflexive Banach
spaces with unconditional bases}

\author{W. B. Johnson}

\address{Department of Mathematics, Texas A\&M University, College Station,
Texas, 77843}

\email{johnson@math.tamu.edu}

\thanks{Both authors were supported in part by NSF/DMS 02-00690 and 05-03688.
This paper forms a portion of the doctoral dissertation of the second
author, which is being prepared at Texas A\&M University under the direction
of the first author.}

\author{Bentuo Zheng}

\address{Department of Mathematics, Texas A\&M University, College Station,
Texas, 77843}

\email{btzheng@math.tamu.edu}



\begin{abstract}
We prove that the dual or any quotient of a separable reflexive Banach space with the unconditional tree property has the unconditional tree property. This is used to prove that a separable reflexive Banach space with the unconditional tree property embeds into a reflexive Banach space with an unconditional basis. This solves several long standing open
problems.  In particular, it yields that a quotient of a reflexive Banach space with an unconditional finite dimensional decomposition embeds into a reflexive Banach space with an unconditional basis.
\end{abstract}

\maketitle

\section{Introduction}
It has long been known that Banach spaces with unconditional bases as well as their
subspaces are much better behaved than general Banach spaces, and that many of the
reflexive spaces (including $L_p(0,1), 1<p<\infty$) that arise naturally in analysis
have unconditional bases. It is however difficult to determine whether a given
Banach space has an unconditional basis or embeds into a space which has an unconditional
basis. Two problems, considered important since at least the 1970's, stand out.
\begin{itemize}
\item[(a)] Give an intrinsic condition on a Banach space $X$ which is equivalent to the
embeddability of $X$ into a space with an unconditional basis.
\item[(b)] Does every complemented subspace of a space with an unconditional basis have an
unconditional basis?
\end{itemize}
Problem (b) remains open, but in this paper we provide a solution to problem (a) for reflexive
Banach spaces. This characterization also yields that a quotient of a reflexive space with an
unconditional basis embeds into a reflexive space with unconditional basis, which solves
another problem from the 1970's. Here some condition on the space with an unconditional basis
is needed because every separable Banach space is a quotient of $\ell_1$.

There is, of course, quite a lot known around problems (a) and
(b). For example, Pe\l czy\'nski and Wojtaszczyk \cite{PW} proved
that if $X$ has an unconditional expansion of identity (i.e., a
sequence $(T_n)$ of finite rank operators such that $\sum T_n$
converges unconditionally in the strong operator topology to the
identity on $X$), then $X$ is isomorphic to a complemented
subspace of a space that has an unconditional finite dimensional
decomposition (UFDD). Later, Lindenstrauss and Tzafriri \cite{LT}
showed that every space with an UFDD embeds (not necessarily
complementably) into a space with an unconditional basis. As
regards reflexive spaces, it was shown in \cite{FJT} using a
result from \cite{DFJP} (and answering a question from that paper)
that if a reflexive Banach space embeds into a space with an
unconditional basis, then it embeds into a reflexive space with an
unconditional basis. As regards the quotient problem we mentioned
above, Feder \cite{F} gave a partial solution by proving that if
$X$ is a quotient of a reflexive space which has an UFDD and $X$
has the approximation property, then $X$ embeds into a space with
an unconditional basis.

It is well known and easy to see that if a Banach space $X$ embeds into a space with an unconditional
basis, then $X$ has the unconditional subsequence property; that is, there exists a $K>0$ so that every
normalized weakly null sequence in $X$ has a subsequence which is $K$-unconditional. In fact, failure
of the unconditional subsequence property is the only known criterion for proving that a given reflexive
space does not embed into a space with an unconditional basis. However, in the last section we
construct a Banach space which has the unconditional subsequence property but does not embed into
a Banach space that has an unconditional basis. This is not surprising, given previous examples of
E. Odell and Th. Schlumprecht \cite{OS1}. Moreover, Odell and Schlumprecht have taught us that by
replacing a subsequence property by the corresponding ``branch of a tree'' property, you get a
stronger property that sometimes can be used to give a characterization of spaces that embed into
a space with some kind of structure. The property relevant for us is the unconditional tree
property and Odell and Schlumprecht's beautiful results are   essential tools for us.

 We use standard
Banach space theory terminology, as can be found in \cite{LT}.

\section{Main results}

\begin{definition}
Let $[\mathbf{N}]^{<\omega}$ denote all finite subsets of the positive integers.
By a normalized weakly null tree, we mean a family
$(x_A)_{A \in [\mathbf{N}]^{<\omega}} \subset S_X$
with the property that every sequence $(x_{A\bigcup\{n\}})_{n\in\mathbf{N}}$
is weakly null. Let $A=\{n_1,...,n_m\}$ with $n_1<...<n_m$ and $B=\{j_1,...,j_r\}$
with $j_1<...<j_r$. Then we say $A$ is an initial segment of $B$ if $m\leq r$
and $n_i=j_i$ when $1\leq i\leq m$. The tree order on
$(x_A)_{A \in [\mathbf{N}]^{<\omega}}$ is given by $x_A\leq x_B$ if $A$ is an initial
segment of $B$. A branch of a tree is a maximal linearly ordered subset of the tree
under the tree order. We say $X$ has the $C$-UTP if every normalized weakly null tree
in $X$ has a $C$-unconditional branch for some $C>0$. $X$ has the UTP if $X$ has the
$C$-UTP for some $C>0$.
\end{definition}

\begin{remark}
E. Odell, Th. Schlumprecht and A. Zsak proved in \cite{OSZ} that
if every normalized weakly null tree in $X$ admits a branch which
is unconditional, then $X$ has the $C$-UTP for some $C>0$. A
simpler proof will appear in the forthcoming paper of R. Haydon,
E. Odell and Th. Schlumprecht \cite{HOS}. So there is no ambiguity
when using the term ``UTP".
\end{remark}

Given an FDD $(E_n)$, $(x_n)$ is said to be a block sequence
with respect to $(E_n)$ if there exists a sequence of integers
$0=m_1<m_2<m_3<...$ such that $x_n \in
\bigoplus_{j=m_{n}}^{m_{n+1}-1} E_j, \forall n\in \mathbf{N}$.
$(x_n)$ is said to be a skipped block sequence with respect to $(E_n)$ if
there exists a sequence of increasing integers $0=m_1< m_2<
m_3<...$ such that $m_n+1<m_{n+1}$ and $x_n \in
\bigoplus_{j=m_{n}+1}^{m_{n+1}-1} E_j, \forall n\in \mathbf{N}$.
Let $\delta=(\delta_i)$ be a sequence of positive numbers
decreasing to $0$. We say $(y_n)$ is a $\delta$-skipped block
sequence with respect to $(E_n)$ if there is a skipped block
sequence $(x_n)$ so that $\|y_n-x_n\|<\delta_n\|y_n\|$ for all
$n\in\mathbf{N}$. We say $(F_n)$ is a blocking of $(E_n)$ if there
is a sequence of increasing integers $0=k_0<k_1<...$ so that
$F_n=\oplus_{j=k_{n-1}+1}^{k_n} E_j$.

\begin{definition}
Let $X$ be a Banach space with an FDD $(E_n)$. If there exists a $C>0$ so that every
skipped block sequence with respect to $(E_n)$ is $C$-unconditional, then we say
$(E_n)$ is an unconditional skipped blocked FDD (USB FDD).
\end{definition}

The following is a blocking lemma of W. B. Johnson and M. Zippin (see \cite{JZ2}
or Proposition 1.g.4(a) in \cite{LT}) which will be used later.

\begin{lemma}\label{olemma}
Let $T:X\rightarrow Y$ be a bounded linear operator. Let $(B_n)$ be a shrinking
FDD of $X$ and let $(C_n)$ be an FDD of $Y$. Let $(\delta_n)$ be a sequence of
positive numbers tending to $0$. Then there are blockings $(B_n^{'})$ of $(B_n)$
and $(C_n^{'})$ of $(C_n)$ so that, for every $x\in B_n^{'}$, there is a
$y\in C_{n-1}^{'}\bigoplus C_n^{'}$ such that $\|Tx-y\|\leq \delta_n\|x\|$.
\end{lemma}

The lemma above actually works for any further blockings of $(B_n^{'})$ and
$(C_n^{'})$. To be more precise, we have the following stronger result which
is actually a formal consequence of Lemma~\ref{olemma} as stated.

\begin{lemma}\label{slemma}
Let $T:X\rightarrow Y$ be a bounded linear operator. Let $(B_n)$ be a shrinking
FDD of $X$ and let $(C_n)$ be an FDD of $Y$. Let $(\delta_n)$ be a sequence of
positive numbers tending to $0$. Then there are blockings $(B_n^{'})$ of $(B_n)$
and $(C_n^{'})$ of $(C_n)$ so that, for any further blockings $(\tilde{B_n})$ of
$(B_n^{'})$ with $\tilde{B_n}=\oplus_{i=k_n}^{k_{n+1}-1} B_i^{'}$ and $(\tilde{C_n})$
of $(C_n^{'})$ with $\tilde{C_n}=\oplus_{i=k_n}^{k_{n+1}-1} C_i^{'}$ and for
any $x\in \tilde{B_n}$, there is a $y\in \tilde{C_{n-1}}\oplus \tilde{C_n}$ such
that $\|Tx-y\|\leq \delta_n\|x\|$.
\end{lemma}

\begin{proof}
Let $(\delta_i)$ be a sequence of positive numbers decreasing to
$0$. Let $(\tilde{\delta_i})$ be another sequence of positive
numbers which go to $0$ so fast that
$\sum_{j=i}^{\infty}\tilde{\delta_j}< 2\lambda\delta_i$, where $\lambda$ is the basis constant for $(B_n)$. By
Lemma~\ref{olemma}, we get blockings $(B_n^{'})$ of $(B_n)$ and
$(C_n^{'})$ of $(C_n)$ so that for every $x\in B_n^{'}$, there is
a $y\in C_{n-1}^{'}\oplus C_n^{'}$ such that
$\|Tx-y\|\leq\tilde{\delta_n}\|x\|$. Let
$\tilde{B_n}=\oplus_{i=k_n}^{k_{n+1}-1} B_i^{'}$ and
$\tilde{C_n}=\oplus_{i=k_n}^{k_{n+1}-1} C_i^{'}$ be blockings of
$(B_n^{'})$ and $(C_n^{'})$. Let $x\in \tilde{B_n}$. Then we can
write $x=\sum_{i=k_n}^{k_{n+1}-1} x_i, x_i\in B_i^{'}$. So by our
construction of $(B_n^{'})$ and $(C_n^{'})$, there are $y_i\in
C_{i-1}^{'}\oplus C_i^{'}, k_n\leq i\leq k_{n+1}-1$ so that
$\|Tx_i-y_i\|\leq\tilde{\delta_i}\|x_i\|, k_n\leq i\leq
k_{n+1}-1$. Let $y=\sum_{i=k_n}^{k_{n+1}-1}
y_i\in\tilde{C_{n-1}^{'}}\oplus\tilde{C_n^{'}}$. Then we have
\begin{equation*}
\|Tx-y\|\leq \sum_{i=k_n}^{k_{n+1}-1}
\tilde{\delta_i}\|x_i\|\leq \sum_{i=k_n}^{k_{n+1}-1}
2\lambda\tilde{\delta_i}\delta_i\|x\|\leq \delta_n \|x\|.
\end{equation*}
\end{proof}

The following convenient reformulation of Lemma~\ref{olemma} will also be used
(see \cite{JZ1} and \cite{JZ2} or \cite{OS2}).

\begin{lemma}\label{lemma}
Let $T:X\mapsto Y$ be a bounded linear operator. Let $(B_n)$ be a shrinking FDD for $X$
and $(C_n)$ be an FDD for $Y$. Let $(\delta_i)$ be a sequence of positive numbers
decreasing to $0$. Then there is a blocking $(B_n')$ of $(B_n)$ and a blocking
$(C_n')$ of $(C_n)$ so that for any $x\in B_n'$ and any $m\neq n, n-1$,
\begin{equation*}
\|Q_m(Tx)\|<\delta_{\max\{m,n\}}\|x\|,
\end{equation*}
where $Q_j$ is the canonical projection from $Y$ onto $C_j'$.
\end{lemma}

\begin{remark}
The qualitative content of Lemma~\ref{lemma} is that there are blockings $(B_n')$ of $(B_n)$  and
$(C_n')$ of $(C_n)$ so that
$TB_n'$ is essentially contained in $C_{n-1}'+C_n'$.
\end{remark}

Our first theorem says that the unconditional tree property for reflexive Banach spaces passes
to quotients. It plays a key role in this paper and involves the lemmas above as well as results
and ideas of Odell and Schlumprecht.

Let us explain the sketch of the proof of the special case when
$Y$ is a reflexive space with the UTP and $Y$ has an FDD $(E_n)$,
while $X$ is a quotient of $Y$ which has an FDD $(V_n)$. Since $Y$
has the UTP, by Odell and Schlumprecht's fundamental result
\cite{OS1}, there is a blocking $(F_n)$ of the $(E_n)$ which is an
USB FDD. Then we use the ``killing the overlap" technique of
\cite{J} to get a further blocking $(G_n)$ so that any norm one
vector $y$ is a small perturbation of the sum of a skipped block
sequence $(y_i)$ with respect to $(F_n)$ and $y_i\in G_{i-1}\oplus
G_i$. Let $Q: Y\mapsto X$ be the quotient map. Using
Lemma~\ref{slemma} and passing to a further blocking, without loss
of generality, we assume that $QG_i$ is essentially contained in
$H_{i-1}+H_i$, where $(H_i)$ is the corresponding blocking of
$(V_n)$. Let $(x_A)$ be a normalized weakly null tree in $X$. We
then choose a branch $(x_{A_i})$ so lacunary that $(x_{A_i})$ is a
small perturbation of a block sequence of $(H_n)$ and for each
$i$, there is at least one $H_{k_i}$ between the essential support
of $x_{A_i}$ and $x_{A_{i+1}}$. Let $x=\sum a_ix_{A_i}$ with
$\|x\|=1$. Considering a preimage $y$ of $x$ under the quotient
$Q$ from $Y$ onto $X$ (with $\|y\|=1$), by our construction, we
can essentially write $y$ as the sum of $(y_i)$ where $(y_i)$ is a
skipped block sequence with respect to $(F_n)$. Since $(F_n)$ is
USB, $(y_i)$ is unconditional. By passing to a suitable blocking
$(z_i)$ of $(y_i)$ and using Lemma~\ref{slemma}, it is not hard to
show that $Qz_i$ is essentially equal to $a_ix_{A_i}$. Noticing
that $(z_i)$ is unconditional, we conclude that $(x_{A_i})$
is also unconditional.

For the general case when $X$ and $Y$ do not have an FDD, we have
to embed them into some superspaces with FDD. The difficulty is
that when we decompose a vector in $Y$ as the sum of disjointly
supported vectors in the superspace, we do not know that the
summands are in $Y$. The same problem occurs for vectors in $X$.
This makes the proof rather technical and a lot of computations
appear.

\begin{theorem}\label{theorem2}
Let $X$ be a quotient of a separable reflexive Banach space $Y$ with the UTP. Then $X$ has the UTP.
\end{theorem}

\begin{proof}
By Zippin's result \cite{Z}, $Y$ embeds isometrically into a reflexive space $Z$
with an FDD. A key point in the proof is that Odell and
Schlumprecht proved (Proposition 2.4 in \cite{OS2}) that there is
a further blocking $(G_n)$ of the FDD for $Z$, $\delta=(\delta_i)$
and a $C>0$ so that every $\delta$-skipped block sequence
$(y_i)\subset Y$ with respect to $(G_i)$ is $C$-unconditional. Let $\lambda$ be the basis constant for $(G_n)$.

Since $X$ is separable, we can regard $X$ as a subspace of
$L_{\infty}$. Let $\epsilon>0$. We may assume that
\begin{itemize}
\item[(a)] $\sum_{j>i}\delta_j<\delta_i$,
\item[(b)] $i\delta_i<\delta_{i-1}$,
\item[(c)] $\sum\delta_i<\epsilon$.
\end{itemize}

 Let $Q$ be a quotient map from $Y$ onto $X$, which
can be extended to a norm one map from $Z$ into $L_{\infty}$ and
we still denote it by $Q$. $QZ$, as any separable subspace of
$L_{\infty}$, is contained in some super space isometric to
$C(\Delta)$ with monotone basis $(v_i)$. Here $\Delta$ is the
Cantor set.

Let $(x_A)$ be a normalized weakly null tree in $X$.
Then we let $(E_n)$ and $(F_n)$ be blockings of $(G_i)$ and $(v_i)$ respectively
which satisfy the conclusions of Lemma~\ref{slemma} and Lemma~\ref{lemma}. Using the
``killing the overlap" technique (see Proposition 2.6 in \cite{OS2}),
we can find a further blocking $(\tilde{E_n}=\oplus_{i=l(n)+1}^{l(n+1)}E_i)$ so that
for every $y\in S_Y$, there exists $(y_i)\subset Y$ and integers $(t_i)$ with
$l(i-1)<t_i\leq l(i)$ for all $i$ such that
\begin{itemize}
\item[(I)] $y=\sum y_i$,
\item[(II)] For $i\in\mathbf{N}$, either $\|y_i\|<\delta_i$ or
$\|\sum_{j=t_{i-1}+1}^{t_i-1} P_j y_i-y_i\|<\delta_i\|y_i\|$,
\item[(III)] $\|\sum_{j=t_{i-1}+1}^{t_i-1} P_j y-y_i\|<\delta_i$,
\item[(IV)] $\|P_{t_i}y\|<\delta_i$ for $i\in\mathbf{N}$,
\end{itemize}
where $P_j$ is the canonical projection from $Y$ onto $E_j$.
Let $\tilde{F_n}=\oplus_{i=l(n)+1}^{l(n+1)}F_i$ and let $\tilde{P_j}$ be the canonical
projection from $X$ onto $\tilde{F_j}$. Since $(x_A)$ is a weakly null tree,
we can pick inductively a branch $(x_{A_i})$ and an increasing sequence of integers
$1=k_0<k_1<...$ such that for any $i\in\mathbf{N}$
\begin{itemize}
\item[(i)] $\|\sum_{j=k_{2i-2}}^{k_{2i-1}-1}\tilde{P_j}x_{A_i}-x_{A_i}\|<\delta_i$,
\item[(ii)] $\|\sum_{j=k_{2i-2}}^{k_{2i-1}-1}\tilde{P_j}x_{A_t}\|<\delta_{\max\{i,t\}}$,
for any $t\neq i$.
\end{itemize}
We will prove that $(x_{A_i})$ is unconditional.
Let $x=\sum a_i x_{A_i}, \|x\|=1$.
Let $y\in S_Y$ so that $Q(y)=x$. Then $y$ can be written as $\sum y_j$ where $(y_j)$ satisfy
(I), (II), (III) and (IV). Define $k_{-1}=-1$ and let $z_i=\sum_{j=k_{2i-3}+2}^{k_{2i-1}+1} y_j$.
We will prove that $\|Qz_i-a_ix_{A_i}\|$ is small.

\begin{equation}\label{equation}
\begin{split}
\|Qz_i-a_ix_{A_i}\|
&\leq\|Q(\sum_{j=t_{k_{2i-3}+1}}^{t_{k_{2i-1}+1}} P_jy)
-(\sum_{j=k_{2i-2}}^{k_{2i-1}-1}\tilde{P_j})x\|\\
&+\|z_i-\sum_{j=t_{k_{2i-3}+1}}^{t_{k_{2i-1}+1}} P_jy\|\\
&+\|a_ix_{A_i}-(\sum_{j=k_{2i-2}}^{k_{2i-1}-1}\tilde{P_j})x\|.
\end{split}
\end{equation}
Hence we need to estimate the three terms in the right hand side of the above inequality.
By the construction, for $i>1$, we have
\begin{equation}\label{equation1}
\begin{split}
\|z_i-\sum_{j=t_{k_{2i-3}+1}}^{t_{k_{2i-1}+1}} P_jy\|
&<\sum_{j=k_{2i-3}+2}^{k_{2i-1}+1}(\|\sum_{l=t_{j-1}+1}^{t_j-1} P_ly-y_j\|+\|P_{t_{j-1}}y\|) + \|P_{t_{k_{2i-1}+1}}y\|\\
&<\sum_{j=k_{2i-3}+2}^{k_{2i-1}+1}\delta_j+\sum_{j=k_{2i-3}+2}^{k_{2i-1}+2}\delta_{j-1}\\
&<\delta_{k_{2i-3}+1}+\delta_{k_{2i-3}}\\
&<\delta_i.
\end{split}
\end{equation}
By direct calculation, for $i=1$, we have
\begin{equation}\label{equation2}
\|z_1-\sum_{j=1}^{t_{k_1+1}} P_jy\|<2\delta_1.
\end{equation}
This gives an estimate of the second term.
For the third term, we have
\begin{equation}\label{equation3}
\begin{split}
\|a_ix_{A_i}-(\sum_{j=k_{2i-2}}^{k_{2i-1}-1}\tilde{P_j})x\|
&<\|(\sum_{j=k_{2i-2}}^{k_{2i-1}-1}\tilde{P_j})(a_ix_{A_i}-x)\|
+\|a_i(x_{A_i}-(\sum_{j=k_{2i-2}}^{k_{2i-1}-1}\tilde{P_j})x_{A_i}\|\\
&<2(k_{2i-2}\delta_{k_{2i-2}}+\sum_{j\geq k_{2i-1}}\delta_j)+2\delta_i\\
&<2(\delta_{k_{2i-2}-1}+\delta_{k_{2i-1}-1})+2\delta_i\\
&<4\delta_i.
\end{split}
\end{equation}
For the first term, let $Q_j$ be the canonical projection from $X$ onto $F_j$ and
let $J_1=[t_{k_{2i-3}+1}, t_{k_{2i-1}+1}], J_2=[l_{k_{2i-2}}+1, l_{k_{2i-1}}]$ and
$J_1'=(t_{k_{2i-3}+1}, t_{k_{2i-1}+1})$.
Then we have
\begin{equation}\label{equation4}
\begin{split}
\|Q(\sum_{j\in J_1}P_jy)-(\sum_{j\in J_2}Q_j)Qy\|
&\leq \|Q(\sum_{j\in J_1}P_jy)-(\sum_{j\in J_1}Q_j)Qy\|
+\|(\sum_{j\in J_1}Q_j)Qy-(\sum_{j\in J_2}Q_j)Qy\|\\
&=\|Q(\sum_{j\in J_1}P_jy)-(\sum_{j\in J_1}Q_j)Qy\|
+\|\sum_{j\in J_1-J_2}Q_j(\sum a_ix_{A_i})\|\\
&<\|Q(\sum_{j\in J_1}P_jy)-(\sum_{j\in J_1}Q_j)Qy\|+4\delta_i\\
&\leq\|(\sum_{j\notin J_1}Q_j)Q(\sum_{j\in J_1}P_jy)\|
+\|(\sum_{j\in J_1}Q_j)Q(\sum_{j\notin J_1}P_j y)\|+4\delta_i\\
&<\|(\sum_{j\notin J_1}Q_j)Q(\sum_{j\in J_1'}P_jy)\|
+\|(\sum_{j\in J_1'}Q_j)Q(\sum_{j\notin J_1}P_j y)\|
+6\delta_i\\
&<2\lambda\delta_i + 2\lambda\delta_i + 6\delta_i \\
&=(4\lambda + 6)\delta_i.
\end{split}
\end{equation}
From Inequality~\ref{equation1}, ~\ref{equation2},
~\ref{equation3} and ~\ref{equation4}, we conclude that
\begin{equation*}
\|Qz_i-a_ix_{A_i}\|<(4\lambda + 12)\delta_i.
\end{equation*}
Let $(\epsilon_i)\subset\{-1,1\}^{\mathbf{N}}$. Let
$I\subset\mathbf{N}$ be the set of indices $i\in\mathbf{N}$ for which $\|y_i\|<\delta_i$ and let $I_i=[k_{2i-3}+2, k_{2i-1}+1]$. So
$z_i=\sum_{j\in I_i}y_j$. Let $z_i'=\sum_{j\in I_i-I}y_j$. It is
easy to verify that $\|z_i-z_i'\|<\delta_i$. Hence
$\|Qz_i'-a_ix_{A_i}\|<(4\lambda + 13)\delta_i$. Now by (II), we know that
$(z_i')$ is a $\delta$-skipped block sequence. Hence, $(z_i')$ is
unconditional. So we have
\begin{equation*}
\begin{split}
\|\sum\epsilon_i a_i x_{A_i}\| &\leq \|Q(\sum\epsilon_i z_i')\|+(4\lambda + 13)\epsilon \\
&\leq C\|\sum z_i'\|+ (4\lambda + 13)\epsilon \\
&<C(\|\sum z_i\|+\sum\delta_i)+ (4\lambda + 13)\epsilon\\
&\leq C+ (C+4\lambda + 13)\epsilon,
\end{split}
\end{equation*}
This shows $(x_{A_i})$ is an unconditional sequence.
\end{proof}

\begin{remark}
If the original space $Y$ has the $1+\epsilon$-UTP for any $\epsilon>0$, then any quotient of $Y$ has the $1+\epsilon$-UTP for any $\epsilon>0$.
\end{remark}

The following is an elementary lemma which will be used later. We omit the standard proof.

\begin{lemma}\label{lemma1}
Let $X$ be a Banach space and $X_1, X_2$ be two closed subspace of $X$.
If $X_1\cap X_2=\{0\}$ and $X_1+X_2$ is closed, then $X$ embeds into $X/X_1\oplus X/X_2$.
\end{lemma}

In \cite{JR}, W. B. Johnson and H. P. Rosenthal proved that any separable Banach
space $X$ admits a subspace $Y$ so that both $Y$ and $X/Y$ have a FDD. The proof
uses    Markuschevich bases. A Markuschevich basis for a separable Banach
space $X$ is a biorthogonal system $\{x_n, x_n^*\}_{n\in\mathbf{N}}$ for which
the span of the $x_n$'s is dense in $X$ and the $x_n^*$'s separate the
points of $X$. By Theorem 1.f.4 in \cite{LT}, every separable Banach space $X$ has a
Markuschevich basis $\{x_n, x_n^*\}_{n\in\mathbf{N}}$  so that $[x_n^*]$ contains any
designated separable subspace of $X^*$. The following lemma is a stronger form of the result
of Johnson and Rosenthal which follows from the original proof. For the convenience of
the readers, we give a sketch of the proof. We use $[x_i]_{i\in I}$ to denote the closed
linear span of $(x_i)_{i\in I}$.

\begin{lemma}\label{lemma2}
Let $X$ be a separable Banach space. Then there exists a subspace $Y$ with FDD $(E_n)$
so that for any blocking $(F_n)$ of $(E_n)$ and for any sequence $(n_k)\subset\mathbf{N}$,
$X/\overline{\textrm{span}\{(F_{n_k})_{k=1}^{\infty}\}}$ admits an FDD $(G_n)$. Moreover,
if $X^*$ is separable, $(E_n)$ and $(G_n)$ can be chosen to be shrinking.
\end{lemma}

\begin{proof}
Let $\{x_i, x_i^*\}$ be a Markuschevich basis for $X$ so that  $[x_i^*]$ is a norm
determining subspace of $X^*$ and even
$[x_i^*]=X^*$ if
$X^*$ is separable. Then we can choose inductively finite sets
$\sigma_1\subset\sigma_2\subset...$ and
$\eta_1\subset\eta_2\subset...$ so that $\sigma=\bigcup_{n=1}^{\infty}\sigma_n$ and
$\eta=\bigcup_{n=1}^{\infty}\eta_n$ are complementary infinite subsets of the positive
integers and for $n=1, 2,...$,
\begin{itemize}
\item[(i)] if $x^*\in [x_i^*]_{i\in\eta_n}$, there is a $x\in [x_i]_{i\in\eta_n\bigcup\sigma_{n+1}}$
so that $\|x\|=1$ and $|x^*(x)|>(1-\frac{1}{n+1})\|x^*\|$;
\item[(ii)] if $x\in [x_i]_{i\in\sigma_n}$, there is a $x^*\in [x_i^*]_{i\in\sigma_n\bigcup\eta_n}$
so that $\|x^*\|=1$ and $|x^*(x)|>(1-\frac{1}{n+1})\|x\|$.
\end{itemize}
Once we have this, by the proof of  Theorem IV.4 in \cite{JR}, we have
$[x_i]_{i\in\sigma}^{\perp}$ is the $w^*$ closure of $[x_i^*]_{i\in\eta}$.
Put $Y=[x_i^*]_{i\in\eta}^{\perp}=[x_i]_{i\in\sigma}$. By the analogue of Proposition II.1(a) in
\cite{JR}, we deduce that  $X/Y$    has an FDD and that  $([x_i]_{i\in\sigma_n})_{n=1}^{\infty}$
forms an FDD for $Y$. So to prove Lemma~\ref{lemma2}, it is enough to prove that for any blocking
$(\Sigma_n)$ of $(\sigma_n)$ or any subsequence $(\sigma_{n_k})$ of $(\sigma_n)$ (this of course
needs the redefining of $(\eta_n)$), (i) and (ii) still hold. But this is more or less obvious
because if
$\Sigma_n=\bigcup_{i=k_{n-1}+1}^{k_n}\sigma_i$, then we define $\Delta_n=\bigcup_{i=k_{n-1}+1}^{k_n}\eta_i$
and it is easy to check $\{\Sigma_n, \Delta_n\}$ satisfy (i) and (ii). For a subsequence
$(\sigma_{n_k})$, if we let $\Sigma_k=\sigma_{n_k}$ and define
$\Delta_k=\bigcup_{i=n_k}^{n_{k+1}-1}\eta_i$, then $\{\Sigma_n, \Delta_n\}$
satisfy (i) and (ii). The rest is exactly the same as in the proof of Theorem IV.4 in \cite{JR}.

\end{proof}

The next lemma shows that for a reflexive space with an USB FDD, its dual also has an USB FDD.

\begin{lemma}\label{alemma2}
Let $X$ be a reflexive Banach space with an USB FDD $(E_n)$. Then there is a
blocking $(F_n)$ of $(E_n)$ so that $(F_n^*)$ is an USB FDD for $X^*$.
\end{lemma}

\begin{proof}
Without loss of generality, we assume $(E_n)$ is monotone.
Let $(\delta_i)$ be a sequence of positive numbers deceasing fast to $0$.
By the ``killing the overlap" technique, we get a blocking $(F_n)$ of $(E_n)$
with $F_n=\sum_{i=k_{n-1}+1}^{k_n} E_i$ so that given any $x=\sum x_i$ with
$x_i\in E_i, \|x\|=1$, there is an increasing sequence $(t_n)$  with
$k_{n-1}<t_n<k_n$ such that $\|x_{t_i}\|<\delta_i$, where $0=k_0<k_1<...$.
Let $(F_n^*)$ be the dual FDD of $(F_n)$ and let $(x_i^*)$ be a normalized
skipped block sequence with respect to $(F_n^*)$ so that
$x_i^*\in\oplus_{j=m_{i-1}+1}^{m_i-1}F_j^*$ where $0=m_0<m_1<...$.
Let $x^*=\sum a_i x_i^*$ with $\|x^*\|=1$. Let $x=\sum x_i$ be a norming
functional of $x^*$ with $x_i\in E_i$. By the definition of $(F_n)$, we get
an increasing sequence $(t_i)$ with $k_{i-1}<t_i<k_i$ so that $\|x_{t_i}\|<\delta_i$.
We define $y_1=\sum_{j=1}^{t_{m_1}-1} x_j$ and $y_i=\sum_{j=t_{m_{i-1}}+1}^{t_{m_i}-1} x_j$
for $i>1$. Let $y=\sum y_i$. So by the triangle inequality,
\begin{equation*}
\|x-y\|\leq\|\sum x_{t_{m_i}}\|\leq\sum\|x_{t_{m_i}}\|<\sum\delta_{t_{m_i}}.
\end{equation*}
Let $(\epsilon_i)\subset\{-1, 1\}^{\mathbf{N}}$ and let
$\tilde{x^*}=\sum\epsilon_i a_ix_i^*$. We will estimate
$\tilde{x^*}(\sum\epsilon_i y_i)$.
\begin{equation*}
\begin{split}
|\tilde{x^*}(\sum\epsilon_i y_i)|
&=|\sum\epsilon_i a_ix_i^*(\sum\epsilon_i y_i)|\\
&=|\sum a_ix_i^*(\sum y_i)|\\
&=|x^*(y)|\\
&\geq 1-\sum\delta_{t_{m_i}}.
\end{split}
\end{equation*}
Since $(y_i)$ is a skipped block sequence with respect to $(E_i)$, $(y_i)$ is
unconditional. Hence
\begin{equation*}
\|\sum\epsilon_i y_i\|\leq C\|\sum y_i\|< C(1+\sum\delta_{m_i}),
\end{equation*}
where $C$ is the unconditional constant associated with the USB FDD $(E_n)$.
If we let $\sum\delta_i<\epsilon/2$, then we conclude that
\begin{equation*}
\|\tilde{x^*}\|>(1-\epsilon)/C(1+\epsilon).
\end{equation*}
Therefore, $(x_i^*)$ is unconditional  with unconditional constant less than $(1+3\epsilon)C$ if $\epsilon$
is sufficiently small. Hence $(F_n^*)$ is an USB FDD.

\end{proof}

\begin{theorem}\label{theorem3}
Let $X$ be a separable reflexive Banach space. Then the following are equivalent.
\begin{itemize}
\item[(a)] $X$ has the UTP.
\item[(b)] $X$ embeds into a reflexive Banach space with an USB FDD.
\item[(c)] $X^*$ has the UTP.
\end{itemize}
\end{theorem}

\begin{proof}
It is obvious that (b) implies (a). If we can prove (a) implies (b), and $X$
satisfies (b), then by Lemma~\ref{alemma2}, $X^*$ is a quotient of a reflexive space
with an USB FDD. So by Theorem~\ref{theorem2}, $X^*$ has the UTP. Hence we only
need to show that (a) implies (b).
Let $X_1$ be a subspace of $X$ with an FDD $(E_n)$ given by Lemma~\ref{lemma2}.
By Proposition 2.4 in \cite{OS2}, we get a blocking $(F_n)$ of $(E_n)$ so that $(F_n)$
is an USB FDD. Let $Y_1=[ F_{4n}]$ and $Y_2=[ F_{4n+2}]$. Then $(F_{4n})$ and
$(F_{4n+2})$ form unconditional FDDs for $Y_1$ and $Y_2$. By Lemma~\ref{lemma2},
$X/Y_i$ has an FDD. Since $X$ has the UTP, by Theorem~\ref{theorem2}, $X/Y_i$ has the
UTP. Now using Proposition 2.4 in \cite{OS2} again, we know that $X/Y_i$ has an USB FDD.
Noticing that $Y_1\cap Y_2=\{0\}$ and $Y_1+Y_2$ is closed, by Lemma~\ref{lemma1},
we have that $X$ embeds into $X/Y_1\oplus X/Y_2$. Hence $X$ embeds into a reflexive
space with an USB FDD.
\end{proof}

\begin{corollary}\label{corollary1}
Let $X$ be a separable reflexive Banach space with the UTP.
Then $X$ embeds into a reflexive Banach space with an unconditional basis.
\end{corollary}

\begin{proof}
By Theorem~\ref{theorem3}, $X$ embeds into a reflexive space $Y$ with an USB FDD $(E_n)$. 
We prove that $Y$ embeds into a reflexive space with an unconditional FDD. Then as was mentioned in the
introduction, $Y$ embeds into a reflexive space with an unconditional basis and so $X$ does.

By Lemma~\ref{alemma2}, there is a blocking $(F_n)$ of $(E_n)$ so that $(F_n^*)$ is an USB 
FDD for $Y^*$. Now let $Y_1=\oplus F_{4n}$ and let $Y_2=\oplus F_{4n+2}$. Then we have $Y_1\cap Y_2=\{0\}$
and $Y_1+ Y_2$ is closed because $(F_{2n})$, being a skipped blocking of $(E_n)$, is unconditional. By
Lemma~\ref{lemma1}, $Y$ embeds into $Y/Y_1\oplus Y/Y_2$. Since $(Y/Y_i)^*$ is isomorphic to $Y_i^{\perp}$,
it is enough to prove $Y_i^{\perp}$ has an unconditional FDD. Let $G_n^*=F_{4n-3}^*\oplus F_{4n-2}^*\oplus
F_{4n-1}^*$. It is easy to see that $(G_n^*)$ forms an FDD for $Y_1^{\perp}$. Noticing that $(G_n)$ is a
skipped blocking of $(F_n^*)$, we conclude that $(G_n)$ is unconditional. Similarly, we can show that
$Y_2^{\perp}$ admits an unconditional FDD. This finishes the proof.
\end{proof}

\begin{corollary}\label{corollary2}
Let $X$ be a quotient of a reflexive Banach space with an unconditional FDD.
Then $X$ embeds into a reflexive Banach space with an unconditional basis.
\end{corollary}

\begin{proof}
Combine Theorem~\ref{theorem2} and Corollary~\ref{corollary1}.
\end{proof}

We mention again that in 1974  Davis,   Figiel,  Johnson and   Pe\l
czy\'nski proved
\cite{DFJP} that a reflexive Banach space $X$ which embeds into a Banach space with
a shrinking unconditional basis embeds into a reflexive space $X$ with an
unconditional basis. The next year,   Figiel,   Johnson and   Tzafriri
\cite{FJT} got a stronger result by removing the shrinkingness of the unconditional
basis in the hypothesis. Our next corollary gives a parallel result for quotients.

\begin{corollary}\label{corollary3}
Let $X$ be a separable reflexive Banach space. If $X$ is a quotient of a Banach space
with a shrinking unconditional basis, then $X$ is isomorphic to a quotient of a reflexive
Banach space with an unconditional basis.
\end{corollary}

\begin{proof}
Since $X$ is a quotient of a Banach space with a shrinking unconditional basis, $X^*$ is
a subspace of a Banach space with an unconditional basis. Hence, by \cite{FJT}, $X^*$ is
isomorphic to a subspace of a reflexive Banach space with an unconditional basis. Therefore,
$X$ is isomorphic to a quotient of a reflexive Banach space with an unconditional basis.
\end{proof}

\begin{remark}
Corollary \ref{corollary3} is different from the result of Figiel, Johnson and Tzafriri
in that the shrinkingness in our result cannot be removed. The reason is more or less
obvious since every separable Banach space is a quotient of $\ell_1$ which has an
unconditional basis.
\end{remark}

Gluing Theorem~\ref{theorem3}, Corollary~\ref{corollary1}, Corollary~\ref{corollary2}
and Corollary~\ref{corollary3} together, we have the following long list of equivalences.

\begin{theorem}\label{theorem4}
Let $X$ be a  separable reflexive Banach space. Then the following are equivalent.
\begin{itemize}
\item[(a)] $X$ has the UTP.
\item[(b)] $X$ is isomorphic to a subspace of a Banach space with an unconditional basis.
\item[(c)] $X$ is isomorphic to a subspace of a reflexive space with an unconditional basis.
\item[(d)] $X$ is isomorphic to a quotient of a Banach space with a shrinking unconditional basis.
\item[(e)] $X$ is isomorphic to a quotient of a reflexive space with an unconditional basis.
\item[(f)] $X$ is isomorphic to a subspace of a quotient of a reflexive space with an
unconditional basis.
\item[(g)] $X$ is isomorphic to a subspace of a reflexive quotient of a Banach space with a
shrinking unconditional basis.
\item[(h)] $X$ is isomorphic to a quotient of a subspace of a reflexive space with an
unconditional basis.
\item[(i)] $X$ is isomorphic to a quotient of a reflexive subspace of a Banach space with a
shrinking unconditional basis.
\end{itemize}
\end{theorem}

\section{Example}

In this section we give an example of a reflexive Banach space for which there exists a
$C>0$ so that every normalized weakly null sequence admits an $C$-unconditional
subsequence while for any $C>0$ there is a normalized weakly null tree such that every
branch is not $C$-unconditional. The construction is an analogue of Odell and Schlumprecht's
example (see Example 4.2 in \cite{OS1}).

We will first construct an infinite sequence of reflexive Banach spaces ${X_n}$. Each
$X_n$ is infinite dimensional and has the property that for $\epsilon>0$, every normalized weakly null 
sequence has a $1+\epsilon$-unconditional basic subsequence, while there is a normalized weakly null tree
for which every branch is at least $C_n$-unconditional and $C_n$ goes to infinity when
$n$ goes to infinity. Then the $\ell_2$ sum of
 $X_n$'s is a reflexive Banach space with the desired property.

Let $[\mathbf{N}]^{\leq n}$ be the set of all subsets of the positive integers with
cardinality less than or equal to $n$. Let $c_{00}([\mathbf{N}]^{\leq n})$ be the space
of sequences with finite support   indexed by $[\mathbf{N}]^{\leq n}$ and denote
its canonical basis by $(e_A)_{A\in [\mathbf{N}]^{\leq n}}$. Let $(h_i)$ be any 
normalized conditional basic sequence which satisfies a block lower $\ell_2$ estimate, 
for example, the boundedly complete basis of
James' space (see problem 6.41 in \cite{FHHSPZ}). Let $\sum a_A e_A$ be an element of
$c_{00}([\mathbf{N}]^{\leq n})$.
Let $(\beta_k)_{k=1}^m$ be disjoint segments.
By a segment in $[\mathbf{N}]^{\leq n}$, we mean a sequence
$(A_i)_{i=1}^{k} \in [\mathbf{N}]^{\leq n}$
with $A_1 = \{n_1 , n_2 , ... , n_l\}, A_2 = \{n_1 , n_2 , ... , n_l , n_{l+1}\},
..., A_k = \{n_1 , n_2 , ... , n_l , ... , n_{l+k-1}\}$, for
some $n_1 < n_2 < ... <n_{l+k-1}$. Let $\beta_k=\{A_{1,k}, A_{2,k},..., A_{j_k,k}\}$
with $A_{i,k}<A_{i+1,k}$ under the tree order in $[\mathbf{N}]^{\leq n}$.
 Now we define $X_n$ to be the completion of
$c_{00}([\mathbf{N}]^{\leq n})$ under the norm

\begin{equation*}
\|\sum a_A e_A\|_{X_n}=
\sup\{(\sum_{k=1}^m(\|\sum_{A_{i,k}\in \beta_k} a_{A_{i,k}} h_i\|)^2)^{1/2}:
(\beta_k)_{k=1}^m\quad\textrm{are disjoint segments}\}.
\end{equation*}

Let $X=(\sum X_n)_2$. We first show that for any $C>0$ there is a normalized weakly null tree
in $X$ so that every branch is at least $C$-unconditional. Let $M$ be so big that the
unconditional constant of $(h_i)_{i=1}^M$ is greater than $C$. Actually the normalized weakly
null tree $(e_A)_{A\in [\mathbf{N}]^{\leq M}}$ in $X_M$ has the property above because any
branch of it is $1$-equivalent to $(h_i)_{i=1}^M$ since $(h_i)$ has a block lower $\ell_2$ estimate. 
Next we need to verify that for every
$\epsilon>0$, every normalized weakly null sequence in $X$ has an $1+\epsilon$-unconditional
basic subsequence. Actually, we will prove that there is a subsequence which is
$1+\epsilon$-equivalent to the unit vector basis of $\ell_2$. By a gliding-hump argument, it is
not hard to verify the following fact.\\

{\bf Fact.} Let $(Y_k)$ be a sequence of reflexive Banach spaces. And let $Y=(\sum
Y_k)_{\ell_2}$. If for every $\epsilon>0, k\in\mathbf{N}$, every normalized weakly null
sequence in $Y_k$ has a subsequence which is $1+\epsilon$-equivalent to the unit vector basis
of $\ell_2$, then for every $\epsilon>0$, every normalized weakly null sequence in $Y$ has a
subsequence which is $1+\epsilon$-equivalent to the unit vector basis of $\ell_2$.\\

Considering the fact, it is enough to show that for every $\epsilon>0, k\in\mathbf{N}$, every
normalized weakly null sequence in $X_k$ has a subsequence which is $1+\epsilon$ equivalent to
the unit vector basis of $\ell_2$. We prove this by induction.

For $k=1$, $X_1$ is isometric to $\ell_2$, so the conclusion is obvious.

Assume the conclusion is true for $X_k$. By the definition of $X_{k+1}$,
$X_{k+1}$ is isometric to $(\sum(\mathbf{R}\oplus X_k))_{\ell_2}$ (where $\mathbf{R}\oplus X_k$
has some norm so that $\{0\}\oplus X_k$ is isometric to $X_k$).
Hence  by hypothesis and the fact we mentioned above, it is easy to see the conclusion is true in
$X_{k+1}$. This finishes the proof.

\begin{remark}
The proof of the corresponding induction step in Example 4.2 in \cite{OS1} is more complicated than
the very simple induction argument in the previous paragraph.  Schlumprecht realized after
\cite{OS1} was written that the induction could be done so simply and his argument works in our
context.
\end{remark}

{\bf Acknowledgements.} The authors thank the referee for his or her corrections, especially for 
their pointing out the imprecision in the initial construction of the example in Section 3.

\end{document}